# Markov and Lagrange spectra for Laurent series in $1/T$ with rational coefficients*

Nikoleta Kalaydzhieva†
*University College London*

March 25, 2020

**Abstract**

The field of formal Laurent series is a natural analogue of the real numbers, and mathematicians have been translating well-known results about rational approximations to that setting. In the framework of power series over the rational numbers, we define and study the Lagrange spectrum, related to Diophantine approximation of irrationals, and the Markov spectrum, related to representation by indefinite binary quadratic forms. We compute both spectra explicitly, and show that they coincide and exhibit no gaps, contrary to what happens over the reals.

## Introduction

A real irrational number $\alpha$, has its "best rational approximation" given by the convergents obtained by truncating its continued fraction expression. In c.1840 Dirichlet showed that there exist infinitely many such good approximations, with respect to denominator. Namely, there exist infinitely many pairs of integers $p, q$, with $q \neq 0$ satisfying

$$\left|\alpha - \frac{p}{q}\right| < \frac{1}{q^2}.$$

Actually the bound can be improved. In 1891 Hurwitz showed that there are infinitely many rational numbers $p/q$, satisfying

$$\left|\alpha - \frac{p}{q}\right| < \frac{1}{\sqrt{5}q^2}.$$

Furthermore, $\sqrt{5}$ is the largest constant that works for all real irrational numbers, meaning if we increase the constant in the denominator further, the statement no longer holds for $\alpha = (1+\sqrt{5})/2$. However, if we exclude $\sqrt{5}$ (and numbers "equivalent to it") we can reduce the gap further to $1/\sqrt{8}q^2$. For $\alpha \in \mathbb{R}/\mathbb{Q}$, we define the *Lagrange constant*, $L(\alpha)$, to be the largest number $L$ such that the inequality

$$\left|\alpha - \frac{p}{q}\right| < \frac{1}{Lq^2},$$

is satisfied by infinitely many rational numbers $p/q$. Running through all real irrationals, we obtain the *Lagrange spectrum*:

$$\mathcal{L} = \{L(\alpha) \mid \alpha \in \mathbb{R}/\mathbb{Q}\}.$$

*This research was supported by ERC grant n° 670239
†Email: zcahndk@ucl.ac.uk

Alternatively, given a doubly infinite sequence of positive integers, say $A = \ldots, a_{-1},\ a_0,\ a_1, \ldots$, we can define the Lagrange constant as the limsup, for $n$ ranging over the integers, of

$$\lambda_n(A) = [a_{n+1},\ a_{n+2}, \ldots] + [0,\ a_n,\ a_{n-1}, \ldots],$$

where $[a_0,\ a_1, \ldots]$ denotes the continued fraction with partial quotients $a_0,\ a_1, \ldots$. Running through all such sequences $A$, we obtain a second definition of the Lagrange spectrum. Interestingly, if we just consider the suprema of $\lambda_n$ for all integers $n$, the set

$$\mathcal{M} = \{\sup_{n \in \mathbb{Z}} \lambda_n(A) \mid A \text{ doubly infinite sequence of positive integers}\}$$

gives the Markov spectrum. For the classical definition, consider a binary quadratic form with real coefficients, $q = ax^2 + bxy + cy^2$ of positive discriminant $d(q) = b^2 - 4ac$. Let $M(q) = \inf|q(x,y)|$ for $x,$ and $y$ taking integer values, not both zero. The *Markov spectrum* is obtained by normalising these minima by the square root of the discriminant and running through all indefinite forms:

$$\mathcal{M} = \left\{ \frac{\sqrt{d(q)}}{M(q)} \;\middle|\; q \text{ a real binary quadratic form with positive discriminant} \right\}.$$

We should remark, that only the finite positive values of the spectrum are taken into account. In particular we exclude real binary quadratic forms that realise 0 for integers $x, y$ not both zero. This was studied by Markov in [5] and in particular, he showed that for elements below 3, the Markov and Lagrange spectra coincide. For the numbers in the Markov spectrum, greater than 3, a lot less is known. However, we do know that there are intervals which contain no point of $\mathcal{M}$. In particular, the Markov spectrum has gaps and contains but is not equal to the Lagrange spectrum. An extensive survey of the results is given by Cusick and Flahive in [2].

In this paper we work in the setting of formal Laurent series in $1/T$ with rational coefficients. We explicitly compute the Lagrange and Markov Spectra. Furthermore, we prove that the two spectra are identical and that they do not exhibit gaps, i.e

**Theorem.** *The Lagrange Spectrum for* $\mathbb{Q}((1/T))$ *is equivalent to the Markov spectrum for* $\mathbb{Q}((1/T))$*, and is equal to* $\mathbb{N} \cup \{\infty\}$.

The statement of the theorem is a combination of corollary 1, 3 and 4.

A detailed survey of results on Diophantine approximation in fields of power series, some of which we recall in section 1, is given by Lasjaunias in [4]. The article, however concentrates on the *approximation spectrum* of an irrational element of $k((1/T))$, for a finite field $k$, defined by Schmidt in [8]. Its upper bound, called the *approximation exponent* $r(\alpha)$ is such that, given $\varepsilon > 0$

$$\left|\alpha - \frac{P}{Q}\right| < |Q|^{-r(\alpha)-1+\varepsilon}$$

is satisfied by infinitely many rationals $P/Q$, but only finitely many satisfy

$$\left|\alpha - \frac{P}{Q}\right| < |Q|^{-r(\alpha)-1-\varepsilon}.$$

In other words it measures the quality of the approximation of $\alpha$ in terms of the exponent.

Recently, some work has been done on the Lagrange spectrum in the setting of formal Laurent series over finite fields, by Parkkonen and Paulin, and Bugeaud in [7] and [1], respectively. They define and study the nonarchimedian quadratic Lagrange spectrum, whose elements are approximations by the orbit of a given quadratic irrational in $\mathbb{F}_q((T^{-1}))$. In particular, they give analogies to the well-known results over the reals about the closedness and boundedness of the spectrum, as well as computations of its maximum.

Organisation of the paper: in section 1 we set up the scene for the Lagrange spectrum over our setting of formal Laurent series in $1/T$ with rational coefficients and defines the equivalent notions of the continued fractions algorithm, convergents and rational approximations. In section 2 we compute explicitly the Lagrange constant for several sets of examples of quadratic irrationals of even degreed polynomials and describe the Lagrange spectrum. In section 3 we develop the theory of indefinite binary quadratic forms in the setting of formal Laurent series in $1/T$ with rational coefficients and show that analogous results, to the ones of real indefinite binary quadratic forms, hold. In section 4 we prove results on the representation of formal Laurent series by indefinite binary quadratic forms and give a function field equivalent to the classical definition of the Markov spectrum. The paper concludes with section 5 by showing that, in this setting we also have an alternative description of the spectra, via doubly infinite sequences. Furthermore, we use these different forms to show that the Lagrange and Markov spectra coincide.

The results in section 3 and the theorems in section 4 regarding the representation of Laurent series by indefinite quadratic forms, follow an analogous approach to Dickson[3], however there are essential differences in the details.

The author would like to thank James Cann for bringing these problems to her attention and to Ardavan Afshar for the helpful discussions.

# 1 Continued fractions for Laurent series with rational coefficients

Let $\mathbb{Q}[T]$ be the ring of polynomials with coefficients in the rationals, and $\mathbb{Q}(T) = \{A/B \mid A, B \in \mathbb{Q}[T],\ B \neq 0\}$ be its field of fractions. Furthermore,

$$\mathbb{Q}\left((1/T)\right) = \left\{ \sum_{i=-\infty}^{m} a_i T^i \mid m \in \mathbb{Z}, a_i \in \mathbb{Q}, \forall i,\ a_m \neq 0 \right\}$$

will denote the set of formal Laurent series in $1/T$ with coefficients in the rationals.

We can extend the usual definition of degree to $\mathbb{Q}\left((1/T)\right)$ in the following way

**Definition 1.** For $\alpha = \sum_{-\infty}^{m} a_i T^i$, $a_m \neq 0$, define

$$\begin{aligned} \partial eg\ &: \mathbb{Q}\left((1/T)\right) \mapsto \mathbb{Z} \\ &\alpha \mapsto m. \end{aligned}$$

Furthermore, we have the convention that $\partial eg\ 0 = -\infty$.

This map is well defined on rational functions and it agrees with the usual definition of degree, on polynomials, i.e

**Lemma 1.** *For $A, B \in \mathbb{Q}[T]$, with $B \neq 0$, of degrees $m, n$, respectively*

1. $\partial eg\ \frac{A}{B} = \deg A - \deg B$.
2. $\partial eg\ A = \deg A$.

*Proof.* Observe that the first implies the second, so it suffices to prove 1. Consider $A, B \in \mathbb{Q}[T]$, with $B \neq 0$.

$$\begin{aligned} A &= \sum_{i=0}^{m} a_i T^i = a_m T^m \left(1 + \sum_{i=0}^{m-1} A_i T^{i-m}\right) \\ B &= \sum_{i=0}^{n} b_i T^i = b_n T^n \left(1 + \sum_{i=0}^{n-1} B_i T^{i-n}\right) \end{aligned}$$

and $a_m, b_n \neq 0$. Furthermore,

$$\begin{aligned} \frac{A}{B} &= \frac{a_m}{b_n} T^{m-n} \left(1 + O(T^{-1})\right) \left(1 + O(T^{-1})\right)^{-1} \\ &= \frac{a_m}{b_n} T^{m-n} \left(1 + O(T^{-1})\right). \end{aligned}$$

Hence, $\partial eg\ \frac{A}{B} = m - n = \deg A - \deg B$, as required. □

**Remark 1.** *For $\alpha \in \mathbb{Q}((1/T))$, we have that $ord(\alpha) := -\partial eg\ \alpha$ is a valuation. Furthermore, $\mathbb{Q}((1/T))$ is the completion of $\mathbb{Q}(T)$ under it.*

Before we describe the continued fractions algorithm over $\mathbb{Q}((1/T))$, we need to make a final definition

**Definition 2.** The *polynomial part* of $\alpha = \sum_{-\infty}^{m} a_i T^i \in \mathbb{Q}((1/T))$ is given by

$$\lfloor \alpha \rfloor := \begin{cases} 0, & \text{if } \partial eg\ \alpha < 0 \\ \sum_{i=0}^{m} a_i T_i, & \text{if } \partial eg\ \alpha = m > 0. \end{cases}$$

The *fractional part* of $\alpha \in \mathbb{Q}((1/T))$ is defined as $\{\alpha\} := \alpha - \lfloor \alpha \rfloor$.

## 1.1 Continued fraction algorithm over $\mathbb{Q}((1/T))$

Let $\alpha \in \mathbb{Q}((1/T))$. The continued fraction algorithm over function fields works in a similar fashion to the one over the reals. First set $\alpha_0(T) = \alpha(T) \in \mathbb{Q}((1/T))$. Then we define $a_0(T) := \lfloor \alpha_0(T) \rfloor$. Hence $\alpha_0(T) = a_0(T) + \{\alpha_0(T)\}$, with $\{\alpha_0(T)\} \in \mathbb{Q}((1/T))$ of finite negative degree. Therefore $\{\alpha_0(T)\}^{-1}$, also an element of $\mathbb{Q}((1/T))$, is well defined and of positive degree. Set $\alpha_1(T) := \{\alpha_0(T)\}^{-1}$, then $\alpha_0 = a_0 + 1/\alpha_1$. We proceed by recursion. Define

$$\begin{aligned} a_i(T) &: = \lfloor \alpha_i(T) \rfloor \\ \alpha_{i+1}(T) &: = \{\alpha_i(T)\}^{-1} \\ \Rightarrow \alpha_i &= a_i + \frac{1}{\alpha_{i+1}}. \end{aligned}$$

Hence

$$\alpha = a_0 + \cfrac{1}{a_1 + \cfrac{1}{a_2 + \ddots \cfrac{1}{\alpha_{i+1}}}},$$

or equivalently

$$\alpha = [a_0,\ a_1,\ \ldots,\ a_i,\ \alpha_{i+1}].$$

The algorithm terminates if the fractional part $\{\alpha_i(T)\}$ is ever 0.

The rational polynomials $a_i$ are called *the partial quotients* of $\alpha$.

**Remark 2.** *The polynomials $a_i(T)$, defined for $i$ up to the point of termination, are all of positive degree, except perhaps for $i = 0$. The partial quotient $a_0(T)$ can be a constant, however the rest must have at least a linear term, since $\partial eg\ a_i(T) = \partial eg\ \lfloor \alpha_i(T) \rfloor = -\partial eg\ \{\alpha_i(T)\} > 0$.*

The continued fraction of $\alpha$ will be infinite for most $\alpha \in \mathbb{Q}((1/T))$. In fact we have

**Proposition 1.** *The continued fraction of $\alpha \in \mathbb{Q}((1/T))$ has finite number of terms, if and only if $\alpha \in \mathbb{Q}(T)$.*

Since the Euclidean algorithm works in $\mathbb{Q}[T]$, the proof of the proposition is identical to the one over the reals.

**Lemma 2.** *Let $D \in \mathbb{Q}[T]$ be a non-square, monic polynomial of even degree. Then $D$ is a square in $\mathbb{Q}((1/T))$, i.e $\sqrt{D}$ has a Laurent series expansion in $1/T$ with rational coefficients.*

*Proof.* Suppose we have $D$ as above, then it must be of the form

$$\begin{aligned} D(T) &= T^{2d} + \sum_{i=0}^{2d-1} a_i T^i \\ &= T^{2d}\left(1 + \sum_{i=0}^{2d-1} a_i T^{i-2d}\right). \end{aligned}$$

Since $\partial eg \left(\sum_{i=0}^{2d-1} a_i T^{i-2d}\right) < 0$, then

$$\left(1 + \sum_{i=0}^{2d-1} a_i T^{i-2d}\right)^{1/2} = \sum_{n=0}^{\infty} \binom{1/2}{n} \left(\sum_{i=0}^{2d-1} a_i T^{i-2d}\right)^n$$

converges in $\mathbb{Q}((1/T))$. Thus $\sqrt{D} = T^d \left(1 + \sum_{i=0}^{2d-1} a_i T^{i-2d}\right)^{1/2}$ is indeed an element of $\mathbb{Q}((1/T))$. □

**Remark 3.** *Notice that we don't necessarily need $D$ to be monic. As long as the leading coefficient of $D$ is a square in $\mathbb{Q}$ the above lemma still holds.*

**Remark 4.** *Since $\sqrt{D} \in \mathbb{Q}((1/T))$, we can compute its continued fraction using the algorithm defined above. Furthermore, since $D$ is not a perfect square, $\sqrt{D} \notin \mathbb{Q}(T)$ and its continued fraction is infinite. However, unlike in the case over the reals, the continued fraction of $\sqrt{D}$ will not always be periodic.*

## 1.2 Convergents

Given an infinite continued fraction expansion, we can truncate at any point, say $[a_0,\ a_1, \ldots,\ a_h]$, and the resulting expression will be a rational function of the form $p_h/q_h(T)$. Furthermore, we can iterate these by the following matrix identity

$$\begin{pmatrix} a_0 & 1 \\ 1 & 0 \end{pmatrix} \begin{pmatrix} a_1 & 1 \\ 1 & 0 \end{pmatrix} \cdots \begin{pmatrix} a_h & 1 \\ 1 & 0 \end{pmatrix} = \begin{pmatrix} p_h & p_{h-1} \\ q_h & q_{h-1} \end{pmatrix}$$

with $p_{-1} = 1,\ q_{-1} = 0$. This provides a sequence of *continuants* $(p_h)_{h\geq 0}$ and $(q_h)_{h\geq 0}$ and thus *convergents* $p_h/q_h$. This very nice matrix representation was given by Van der Poorten and Shallit in [9].

**Proposition 2.** *Given* $p_h/q_h = [a_0,\ a_1, \ldots,\ a_h]$ *and* $p_{h-1}/q_{h-1} = [a_0,\ a_1, \ldots,\ a_{h-1}]$, *we have*

$$\frac{p_h}{p_{h-1}} = [a_{h-1},\ a_{h-2}, \ldots,\ a_0] \quad \textit{and} \quad \frac{q_h}{q_{h-1}} = [a_{h-1},\ a_{h-2}, \ldots,\ a_1]$$

The proof is a direct computation using the recurrence relations connecting the $p_h$'s and $q_h$'s.

$$\begin{aligned} \frac{p_h}{p_{h-1}} &= a_h + \frac{p_{h-2}}{p_{h-1}} \\ &= a_h + \frac{1}{a_{h-1} + \frac{p_{h-3}}{p_{h-2}}} \\ &\cdots \end{aligned}$$

ending at $a_1 = p_0/p_1$. The same computation works for $q_h$ as the same recurrence relation holds, except the final term will be $a_1 = q_1/q_2$, since $q_0 = 0$.

Since we can write $\alpha = [a_0,\ a_1, \ldots,\ a_h,\ \alpha_{h+1}]$ we have the *convergents correspondence*

$$\alpha \leftrightarrow \begin{pmatrix} a_0 & 1 \\ 1 & 0 \end{pmatrix} \cdots \begin{pmatrix} a_h & 1 \\ 1 & 0 \end{pmatrix} \begin{pmatrix} \alpha_{h+1} & 1 \\ 1 & 0 \end{pmatrix} = \begin{pmatrix} p_h & p_{h-1} \\ q_h & q_{h-1} \end{pmatrix} \begin{pmatrix} \alpha_{h+1} & 1 \\ 1 & 0 \end{pmatrix}$$

$$\leftrightarrow \frac{p_h \alpha_{h+1} + p_{h-1}}{q_h \alpha_{h+1} + q_{h-1}}.$$

Therefore

$$\alpha = \frac{p_h \alpha_{h+1} + p_{h-1}}{q_h \alpha_{h+1} + q_{h-1}}.$$

Furthermore, if we take the determinants of the matrices above, we show

**Proposition 3.** *Given a continued fraction expansion of a formal Laurent series* $\alpha = [a_0,\ a_1, \ldots]$, *its continuants* $p_h$ *and* $q_h$ *satisfy*

$$(-1)^h = p_{h-1} q_h - p_h q_{h-1}$$

**Remark 5.** *Notice that if* $h$ *is even, then* $(p_{h-1}, q_{h-1})$ *gives the unique solution to* $q_h x - p_h y = 1$, *such that* $\partial eg\ x < \partial eg\ p_h$ *and* $\partial eg\ y < \partial eg\ q_h$. *To see this observe*

$$\begin{aligned} q_h x - p_h y &= p_{h-1} q_h - p_h q_{h-1} \\ q_h (x - p_{h-1}) &= p_h (y - q_{h-1}). \end{aligned}$$

*Since* $(p_h, q_h) = 1$, *we must have some rational polynomial* $f$, *such that*

$$\begin{aligned} x - p_{h-1} &= f p_h \\ y - q_{h-1} &= f q_h. \end{aligned}$$

*If* $\partial eg\ p_{h-1} < \partial eg\ p_h$ *and* $\partial eg\ q_{h-1} < \partial eg\ q_h$, *then* $f = 0$ *and* $(p_{h-1}, q_{h-1})$ *gives the unique solution to the Diophantine equation* $q_h x - p_h y = 1$, *such that* $\partial eg\ x < \partial eg\ p_h$ *and* $\partial eg\ y < \partial eg\ q_h$.

**Proposition 4.** *The continuants satisfy* $\partial eg\ p_h < \partial eg\ p_{h+1}$ *and* $\partial eg\ q_h < \partial eg\ q_{h+1}$ *for* $h \geq 0$.

*Proof.* We prove the result by induction. First suppose $a_0 \neq 0$. Then since $\partial eg\ a_i > 0$ for all $i > 0$ we have

$$\begin{aligned} \partial eg\ p_1 &= \partial eg\ (a_1 a_0 + 1) = \partial eg\ a_1 + \partial eg\ a_0 > \partial eg\ a_0 = \partial eg\ p_0 \\ \partial eg\ q_1 &= \partial eg\ (a_1) > 0 = \partial eg\ q_0 \end{aligned}$$

If $a_0 = 0$, then $\partial eg\ p_1 = \partial eg\ 1 = 0 > -\infty = \partial eg\ 0 = \partial eg\ p_0$ and $\partial eg\ q_1 = \partial eg\ a_1 > 0 = \partial eg\ q_0$.

Next, we suppose $\partial eg\ p_{h-1} < \partial eg\ p_h$ and $\partial eg\ q_{h-1} < \partial eg\ q_h$, then

$$\begin{aligned} \partial eg\ p_{h+1} &= \partial eg\ (a_{h+1} p_h + p_{h-1}) = \partial eg\ a_{h+1} + \partial eg\ p_h > \partial eg\ p_h \\ \partial eg\ q_{h+1} &= \partial eg\ (a_{h+1} q_h + q_{h-1}) = \partial eg\ a_{h+1} + \partial eg\ q_h > \partial eg\ q_h \end{aligned}$$

as required. □

All the results up until now are well known and analogous to those over the reals and can be found in [6], for example. However, in the setting of the paper we can be a bit more precise and give an exact expression for the degree of $q_h$.

**Lemma 3.** *For* $\alpha \in \mathbb{Q}((1/T))$ *with continued fraction* $[a_0,\ a_1, \ldots]$, $a_i \neq 0$, *and nth convergent* $p_h/q_h$, $h \geq 1$, *we have that*

$$\deg q_h = \sum_{i=1}^{h} \deg a_i.$$

*Proof.* The proof is by induction on $h$. Since $q_1 = a_1$ and $q_2 = a_1 a_2 + 1$, the statement follows easily for $h = 1, 2$. Then suppose $\deg q_h = \sum_{i=1}^{h} \deg a_i$. Consider the recurrence relation

$$q_{h+1} = q_h a_{h+1} + q_{h-1}$$

Since $\deg q_{h-1} < \deg q_h$, and $\deg a_{h+1} \geq 1$, the result follows. □

**Proposition 5.** *Suppose* $\alpha \in \mathbb{Q}((1/T))$ *has a continued fraction expansion* $[a_0,\ a_1, \cdots]$ *and convergents* $p_h/q_h$. *Then* $\partial eg\ \alpha = \partial eg\ \frac{p_h}{q_h}$, *and in particular*

$$\partial eg\ \alpha = \begin{cases} \partial eg\ a_0, & \text{if } \partial eg\ \alpha \geq 0 \\ -\partial eg\ a_1, & \text{otherwise} \end{cases}$$

*Proof.* We prove the result by induction once again. Suppose $a_0 \neq 0$. Then $\partial eg\ p_0 = \partial eg\ a_0$, and $\partial eg\ q_0 = 0$, hence $\partial eg\ (p_0/q_0) = \partial eg\ a_0$. Furthermore, $\partial eg\ \alpha = \partial eg\ \lfloor \alpha \rfloor = \partial eg\ a_0$.

If $a_0 = 0$, then $\partial eg\ p_1 = 0$ and $\partial eg\ q_1 = \partial eg\ a_1$, hence $\partial eg\ (p_1/q_1) = -\partial eg\ a_1$. Furthermore, observe that since $a_0 = 0$, $\alpha = \{\alpha\}$, and $a_1 = \lfloor 1/\{\alpha\} \rfloor = -\partial eg\ \alpha$.

From the recurrence relation for $h > 1$, $\partial eg\ p_h = \partial eg\ a_h + \partial eg\ p_{h-1}$ and $\partial eg\ q_h = \partial eg\ a_h + \partial eg\ q_{h-1}$, we have

$$\partial eg\ \frac{p_h}{q_h} = \partial eg\ p_h - \partial eg\ q_h = \partial eg\ p_{h-1} - \partial eg\ q_{h-1} = \partial eg\ a_0,$$

And the last equality follows from the induction hypothesis.

□

# 2 Rational approximation and the Lagrange spectrum

Similarly to the case over the reals, the convergents of $\alpha$ provide a very good rational approximation.

**Proposition 6.** *Suppose $\alpha \in \mathbb{Q}((1/T))$ and $p, q \in \mathbb{Q}[T]$, with $q \neq 0$. Then*

$$\partial eg\ \left(\alpha - \frac{p}{q}\right) < -2\deg q$$

*if and only if $p/q$ is a convergent for $\alpha$.*

Notice that $p/q$ is a convergent of $\alpha = [a_0,\ a_1 \ldots]$ if and only if $p/q = [a_0,\ a_1, \ldots,\ a_i]$, for some $i \geq 0$. Then the proposition is a direct corollary of the following

**Proposition 7.** *Suppose we have $\alpha, \beta \in \mathbb{Q}((1/T))$, distinct. Then*

$$\deg(\alpha - \beta) < -2\deg q_i,$$

*where $q_i$ is the denominator of the ith convergent of $\alpha$, if and only if the first $i+1$ partial quotients of their continued fraction expansions are the same.*

*Proof.* Suppose $\alpha = [a_0,\ a_1, \ldots,\ a_i,\ \alpha_{i+1}]$, and $\beta = [a_0,\ a_1, \ldots,\ a_i,\ \beta_{i+1}]$, with $\alpha_{i+1} \neq \beta_{i+1}$. Without loss of generality, we can take $\partial eg\ \alpha_{i+1} \leq \partial eg\ \beta_{i+1}$. Then the first $i$ convergents must be the same for both $\alpha$ and $\beta$. From the convergents correspondence,

$$\alpha = \frac{\alpha_{i+1}p_i + p_{i-1}}{\alpha_{i+1}q_i + q_{i-1}}, \quad \text{and} \quad \beta = \frac{\beta_{i+1}p_i + p_{i-1}}{\beta_{i+1}q_i + q_{i-1}}.$$

Taking the difference and applying proposition 3 yields

$$\alpha - \beta = \frac{(-1)^{i+1}(\alpha_{i+1} - \beta_{i+1})}{(\alpha_{i+1}q_i + q_{i-1})(\beta_{i+1}q_i + q_{i-1})}. \tag{1}$$

Considering the degree of both sides of the equality, and using that $\partial eg\ \alpha_{i+1} = \deg a_{i+1}$ and $\partial eg\ \beta_{i+1} = \deg b_{i+1}$ we get

$$\begin{aligned}\partial eg\ (\alpha - \beta) &= -\left(\deg a_{i+1} + \deg b_{i+1} + 2\deg q_i - \deg(a_{i+1} - b_{i+1})\right)\\ &\leq -\deg a_{i+1} - 2\deg q_i\\ &< -2\deg q_i.\end{aligned}$$

For the inequalities we use that $\deg(a_{i+1} - b_{i+1}) \leq \deg b_{i+1}$, by assumption and $\deg a_{i+1} \geq 1$ by definition. This completes the proof in one direction.

For the converse suppose that $\partial eg\ (\alpha - \beta) < -2 \deg q_i$, and $a_0 = b_0, \ldots, a_{h-1} = b_{h-1}$, but $a_h \neq b_h$ for some $h < i$. Without loss of generality, we will assume that $\deg a_h \leq \deg b_h$. If we do the computation (1) for $h-1$ and consider the degree of both sides of the equality, we get

$$\begin{aligned}\partial eg\ (\alpha - \beta) &= -(\deg a_h + \deg b_h + 2 \deg q_{h-1} - \deg(a_h - b_h)) \\ &< -2 \deg q_i.\end{aligned}$$

After rearranging and applying the result from lemma 3, we have

$$\deg a_h + \deg b_h - \deg(a_h - b_h) > 2 \sum_{j=h}^{i} \deg a_j.$$

Furthermore, by assumption $\deg a_h \leq \deg b_h$, hence

$$2 \deg a_h \geq \deg a_h + \deg b_h - \deg(a_h - b_h).$$

Therefore, $\deg a_h > \sum_{j=h}^{i} \deg a_j$ yielding a contradiction. □

We can actually give an explicit formula for how well, in terms of degree, the convergents of $\alpha$ approximate it.

**Theorem 1.** *Suppose $\alpha \in \mathbb{Q}((1/T))$ and $p_h/q_h$ be its $h^{th}$ convergent. Then*

$$\partial eg\ \left(\alpha - \frac{p_h}{q_h}\right) = -2 \deg q_h - \deg a_{h+1}.$$

*Proof.* Let $p_h/q_h$ be the $h^{\text{th}}$ convergent of $\alpha$, then

$$\alpha - \frac{p_h}{q_h} = \frac{(-1)^{h+1}}{q_h(\alpha_{h+1} q_h + q_{h-1})}$$

Considering degree of both sides and using that $\partial eg\ \alpha_{h+1} = \deg a_{h+1}$, by definition, we get

$$\partial eg\ \left(\alpha - \frac{p_h}{q_h}\right) = -2 \deg q_h - \deg a_{h+1}.$$

□

There is no corresponding result to theorem 1 over the reals. Having this identity significantly simplifies, for example the proof of the analogous result to Dirichlet's rational approximation theorem.

**Proposition 8.** *Given $\alpha \in \mathbb{Q}((1/T))$ not a rational function, there exist infinitely many pairs of rational polynomials $p, q$, with $q \neq 0$ such that*

$$\partial eg\ \left(\alpha - \frac{p}{q}\right) \leq -2 \deg q - 1.$$

*Proof.* Since degree of $a_{n+1}$ is always greater or equal to 1, theorem 1 implies that the convergents $p_n/q_n$ satisfy the inequality. Furthermore, since $\alpha$ is not a rational function, proposition 1 implies that there are infinitely many of those. $\square$

If we consider all non-rational $\alpha \in \mathbb{Q}((1/T))$ then we cannot improve the inequality in proposition 8. However, given a specific non-rational Laurent series $\alpha$, we might be able to sharpen the bound. This leads us to the following definition:

**Definition 3.** Given $\alpha \in \mathbb{Q}((T^{-1}))$, we define the *approximation (Lagrange) constant*, $l(\alpha)$ to be the greatest integer $k$ such that

$$\partial eg\ \left(\alpha - \frac{p}{q}\right) \leq -2\deg q - k$$

is satisfied by infinitely many rational polynomials $p,\ q$. We then define the *Lagrange spectrum* over $\mathbb{Q}((1/T))$ to be

$$\mathbb{L} := \{l(\alpha) : \alpha \in \mathbb{Q}((1/T)),\ \text{not rational}\}.$$

**Remark 6.** *From corollary 6, the inequality $\partial eg\ (\alpha - p/q) < -2\partial eg\ q$ is satisfied only by the convergents of $\alpha$, say $p_h/q_h$. By theorem 1, the left hand-side is simply equal to $-2\partial eg\ q_h - \partial eg\ a_{h+1}$, where $a_{h+1}$ is a partial quotient of $\alpha$. Thus we can substantially simplify the definition of the Lagrange constant to*

$$l(\alpha) = \limsup_{h\to\infty} \partial eg\ a_h.$$

**Example 1.** *Suppose $D \in \mathbb{Q}[T]$ is a quadratic polynomial, not a perfect square. Say, $D = (aT+b)^2 + c$, with $a, b, c \in \mathbb{Q}$ and $ac \neq 0$, then*

$$\sqrt{D} = \sqrt{(aT+b)^2 + c} = \left[aT+b,\ \overline{\frac{2}{c}(aT+b),\ 2(aT+b)}\right].$$

*Notice that all partial quotients have degree 1. Therefore $l(\sqrt{D}) = 1$ for $D$ a square-free quadratic polynomial with rational coefficients.*

For more interesting examples of Lagrange constants we need to find $\alpha \in \mathbb{Q}((1/T))$, such that $\deg a_h = d > 1$, for infinitely many $h$.

**Theorem 2.** *For $a, b, c \in \mathbb{Q}[T]$, we have*

1. $\sqrt{a^2+1} = [a,\ \overline{2a}]$;
2. $\sqrt{a^2+c} = [a,\ \overline{2b,\ 2a}]$, *if $a = bc$.*

*Proof.* Observe that part 1 is a consequence of part 2, if we take $b = a$. Hence it suffices to prove the second result. Suppose we are given the continued fraction expansion $[a,\ \overline{2b,\ 2a}] = \alpha \in \mathbb{Q}((1/T))$. This is equivalent to the expression

$$\alpha = a + \frac{1}{\beta}, \text{ where}$$
$$\beta = 2b + \cfrac{1}{2a + \cfrac{1}{\beta}}$$

After rearranging and simplifying the above, we get the following quadratic equation in $\beta$:

$$a\beta^2 - 2ab\beta - b = 0$$
$$\Rightarrow \beta = \frac{ab + \sqrt{a^2b^2 + ab}}{a}.$$

Therefore

$$\begin{aligned}\alpha &= a + \frac{a}{ab + \sqrt{a^2b^2 + ab}} \times \frac{ab - \sqrt{a^2b^2 + ab}}{ab - \sqrt{a^2b^2 + ab}} \\ &= \sqrt{a^2 + \frac{a}{b}} \\ &= \sqrt{a^2 + c}, \text{ where } a = bc.\end{aligned}$$

□

**Example 2.** *Let $d$ be a positive integer. Then the theorem 2 gives us the following examples*

1. $\sqrt{T^{2d} + T^l} = [T^d, \overline{2T^{d-l}, 2T^d}]$, *for* $0 \le l < d$;
2. $\sqrt{T^{2d} + T^d} = [T^d, \overline{2T^d}]$.

**Theorem 3.** *Let $d$ be a positive integer, then*

1. *for $D = T^{2d} + T^l$, where $0 \le l < d$, the continued fraction expansion of $\sqrt{D}$ has partial quotients, $a_h$ with*
$$\deg a_h = \begin{cases} d, & \text{if } h \text{ even} \\ d - l, & \text{if } h \text{ is odd.} \end{cases}$$
2. *for $D = T^{2d} + T^d$, the continued fraction expansion of $\sqrt{D}$ has partial quotients $a_h$ of degree $d$ for all $h \ge 0$.*

*Furthermore, $l(\sqrt{D}) = d$, for any of the polynomials $D$ described in the statement of the theorem.*

*Proof.* Since $D$ is a rational polynomial of even degree, $\sqrt{D} \in \mathbb{Q}((1/T))$ and thus it has an infinite continued fraction expansion. From part *2* of example 2, we see that $\deg a_h = d$, for all $h$, and part *1* of example 2 gives

$$\deg a_h = \begin{cases} d, & \text{if } h \text{ even} \\ d-l, & \text{if } h \text{ is odd.} \end{cases}$$

Finally, remark 6 says $l(\alpha) = \limsup_{h\to\infty} \partial eg\ a_h$, and since $d - l < d$, we conclude $l(\sqrt{D}) = d$ for both parts. □

**Corollary 1.** *The Lagrange Spectrum of* $\mathbb{Q}((1/T))$ *is equal to* $\mathbb{N} \cup \{-\infty\}$.

*Proof.* For each positive integer $k$, there exists $\alpha \in \mathbb{Q}((T^{-1}))$ such that $l(\alpha) = k$. Just take $\alpha$, to be one of the square roots described in theorem 3. □

# 3 Binary quadratic forms over $\mathbb{Q}((1/T))$

We now proceed to set the scene for the definition of the Markov spectrum over $\mathbb{Q}((1/T))$. In order to do so, we need to firstly develop the theory of indefinite binary quadratic forms but in the setting of formal Laurent series in $T^{-1}$ with rational coefficients.

**Definition 4.** A binary quadratic form over $\mathbb{Q}((1/T))$ is defined to be an expression

$$Q = Q(X,Y) = (A,B,C) := AX^2 + BXY + CY^2,$$

where $A, B, C \in \mathbb{Q}((1/T))$, not all rational functions in $T$. We define the discriminant to be $D = B^2 - 4AC$, which is also an element of $\mathbb{Q}((1/T))$.

**Definition 5.** We call a binary quadratic form $(A,B,C)$ *indefinite*, if the discriminant $D$ is a square in $\mathbb{Q}((1/T))$. From lemma 2, this is precisely when $D$ is a polynomial of even degree and with leading coefficient a rational square.

For an indefinite binary quadratic form $Q(X,Y)$, we have that $X - \omega Y$ is a factor, where $\omega$ is a root of

$$A\omega^2 + B\omega + C = 0.$$

We define the *first* and *second* roots to be respectively

$$f := \frac{\sqrt{D} - B}{2A} \qquad s := \frac{-\sqrt{D} - B}{2A}. \tag{2}$$

Furthermore, assuming $A \neq 0$ and $f, s \notin \mathbb{Q}(T)$, the Laurent series for $f, s$ and $\sqrt{D}$ uniquely determine $A, B, C$. Observe that $f$ and $s$ are both in $\mathbb{Q}(T)$ if and only if $A, B$ and $C$ are all rational functions in $T$ and $D$ is a perfect square.

Suppose we substitute

$$x = \alpha X + \beta Y \qquad y = \gamma X + \delta Y, \tag{3}$$

with $\alpha,\ \beta,\ \gamma,\ \delta \in \mathbb{Q}[T]$ not all 0, in $q(x,y)$. This takes the binary quadratic form $q(x,y)$ to the binary quadratic form $Q(X,Y)$. We can also use the matrix form

$$H = \begin{pmatrix} \alpha & \beta \\ \gamma & \delta \end{pmatrix},$$

with the convention that applying the matrix to a binary quadratic form is the same as applying the linear transformation (3) to it.

**Definition 6.** We say that two forms $q$ and $Q$ are *equivalent* if such matrix $H$ exists and $\det(H) = \pm 1$.

Furthermore, we say that $q$ and $Q$ are *properly equivalent* if $\det H = 1$, and *improperly equivalent* if $\det H = -1$.

**Proposition 9.** *The form $q = (a,b,c)$ is transformed into the form $Q = (A,B,C)$ via $H = \begin{pmatrix} \alpha & \beta \\ \gamma & \delta \end{pmatrix} \in GL_2(\mathbb{Q}[T])$, if and only if their first roots $f$ and $F$ and their second roots $s$ and $S$, respectively, are connected by the relations*

$$f = \frac{\alpha F + \beta}{\gamma F + \delta} \quad \textit{and} \quad s = \frac{\alpha S + \beta}{\gamma S + \delta}$$

Analogous computation to the one over the reals works, see [3].

## 3.1 Reduced indefinite binary quadratic forms

**Definition 7.** The indefinite binary quadratic form $Q = (A,B,C)$ is called *reduced* if

$$\partial eg\ f < 0 < \partial eg\ s, \text{ and } f \neq 0.$$

From (2), this is equivalent to

$$\partial eg\ (\sqrt{D} - B) < \partial eg\ (A) < \partial eg\ (\sqrt{D} + B), \text{ and } \sqrt{D} \neq B.$$

**Proposition 10.** *If $q = (A,B,C)$ is reduced, then so is $Q = (C,B,A)$.*

*Proof.* Consider the transformation $\begin{pmatrix} 0 & 1 \\ 1 & 0 \end{pmatrix}$ taking $q$ to $Q$, and in particular the roots $(f,s)$ to $(F,S) = \left(\frac{1}{s}, \frac{1}{f}\right)$. Since $q$ is reduced, then $\partial eg\ f < 0 < \partial eg\ s$. Hence $\partial eg\ F = -\partial eg\ s < 0$, and $\partial eg\ S = -\partial eg\ f > 0$. □

**Theorem 4.** *An indefinite binary quadratic form is properly equivalent to a reduced one.*

*Proof.* Let $q = (a, b, c) = ax^2 + bxy + cy^2$, with $a, b, c \in \mathbb{Q}((1/T))$, be an indefinite binary quadratic form of discriminant $D \neq 0$. It has first and second root $f = (\sqrt{D} - b)/2a$ and $s = (-\sqrt{D} - b)/2a$, respectively. Firstly we will show that $q$ is either a reduced form or is properly equivalent to a binary quadratic form with a first root of non-negative degree. Suppose the degree of $f$ is negative, then either $q$ is already reduced or $\partial eg\ s \leq 0$. If we are in the latter case, apply the transformation $\begin{pmatrix} 0 & 1 \\ -1 & 0 \end{pmatrix}$. Then, $q$ is properly equivalent to an indefinite binary quadratic form with roots $-1/f$ and $-1/s$, both of positive degree.

Hence $q$ is properly equivalent to a binary quadratic form with roots $(\varphi, \sigma)$, such that $\partial eg\ \varphi \geq 0$. Then we apply the transformation $\begin{pmatrix} 1 & h \\ 0 & 1 \end{pmatrix}$ with $h = \lfloor \varphi \rfloor \in \mathbb{Q}[T]$. This takes the roots $(\varphi, \sigma)$ to $(F, S)$, where $F = \{\varphi\}$ and $S = \sigma - h$. Now if $\lfloor \varphi \rfloor \neq \lfloor \sigma \rfloor$, then $\partial eg\ F < 0$ and $\partial eg\ S > 0$, hence $q$ is properly equivalent to a reduced form.

If $\lfloor \varphi \rfloor = \lfloor \sigma \rfloor$, then consider the continued fraction expansion $\varphi = [a_0,\ a_1 \cdots]$ and $\sigma = [b_0,\ b_1, \cdots]$. Pick the smallest $m$ such that $a_m \neq b_m$, $m > 0$, then we have

$$\varphi = [a_0,\ a_1, \cdots,\ a_{m-1},\ f_m] \quad \text{and} \quad \sigma = [a_0,\ a_1, \cdots,\ a_{m-1},\ s_m]$$

Since $a_m \neq b_m$, then $f_m \neq s_m$, and in particular $\lfloor f_m \rfloor \neq \lfloor s_m \rfloor$. Observe that the convergents for $\varphi$ and $\sigma$ are the same up to and including the $(m-1)^{\text{st}}$ term. Then the transformation $\begin{pmatrix} p_{m-1} & p_{m-2} \\ q_{m-1} & q_{m-2} \end{pmatrix}$ takes $(\varphi, \sigma)$ to $(f_m, s_m)$. Furthermore, this matrix has polynomial entries and is of determinant $(-1)^{m-2}$, i.e 1 or $-1$ depending on the parity of $m$. Apply

$$\begin{pmatrix} (-1)^m & h \\ 0 & 1 \end{pmatrix} \quad \text{with } h = \lfloor f_m \rfloor.$$

This takes $(f_m, s_m)$ to $(F, S)$, where $F = (-1)^m \{f_m\}$ has negative degree and $S = (-1)^m (s_m - h)$ has non-negative degree. Since $\lfloor s_m \rfloor \neq \lfloor f_m \rfloor$, $\partial eg\ S$ is positive, and the new quadratic form is reduced and properly equivalent to $q$. □

The reduction algorithm, given in the proof differs to the one over the reals, however the same theorem still holds.

## 3.2 Chain of reduced forms

All results in this section are the direct analogue to the case over the reals and can be found in [3]. We follow the same approach as Dickson, however the proofs differ in the details.

**Theorem 5.** *Each reduced indefinite binary quadratic form has a unique right neighbouring form.*

*Proof.* Let $Q = (A, B, A_1)$ be an indefinite reduced binary quadratic form of discriminant $D$. The transformation $\Delta = \begin{pmatrix} 0 & 1 \\ -1 & \delta \end{pmatrix}$ takes $Q$ to the equivalent form $Q_1 = (A_1, B_1, A_2)$, such that $B_1 = -B - 2\delta A_1$ and $A_2$ obtained from the discriminant $D$. Furthermore,

$$f \xrightarrow{\Delta} F = \delta - \frac{1}{f}$$
$$s \xrightarrow{\Delta} S = \delta - \frac{1}{s}$$

Since $Q$ is reduced, $\partial eg\ f < 0 < \partial eg\ s$. Take $\delta = \lfloor 1/f \rfloor \in \mathbb{Q}[T]$, which has positive degree. Then $\partial eg\ F = \partial eg\ \{1/f\} < 0$ and $\partial eg\ S = \partial eg\ \delta - 1/s = \partial eg\ \delta > 0$, i.e $Q_1$ is reduced. Observe that if $\delta \neq \lfloor 1/f \rfloor$, then $\partial eg\ F > 0$. Hence $Q_1$ is reduced, only if $\delta$ is chosen to be $\lfloor 1/f \rfloor$. □

**Corollary 2.** *Every reduced form has one and only one reduced left neighbouring form.*

*Proof.* If $(A, B, A_1)$ is reduced, then $(A_1, B, A)$ is reduced as well, by proposition 10. From the theorem above, there is a unique reduced right neighbouring form $(A, B_1, A_2)$. Then by proposition 10 $(A_2, B_1, A)$ is also reduced. Furthermore, it has $(A, B, A_1)$ as its unique right neighbouring form. □

Therefore given an indefinite binary quadratic form of discriminant $D \neq 0$ we can construct a chain of equivalent reduced indefinite binary quadratic forms of the same discriminant, say

$$\cdots, \Phi_{-1}, \Phi_0, \Phi_1, \cdots,$$

where $\Phi_i = \big((-1)^i A_i,\ B_i,\ (-1)^{1+1} A_{i+1}\big)$. The transformation $\Delta_i = \begin{pmatrix} 0 & 1 \\ -1 & \delta_i \end{pmatrix}$ takes $\Phi_i$ to $\Phi_{i+1}$. Furthermore, we have the relation $B_i + B_{i+1} = 2g_i A_{i+1}$, where $g_i = (-1)^i \delta_i$.

Let $f_i = \frac{\sqrt{D} - B_i}{(-1)^i 2A_i}$ and $s_i = \frac{\sqrt{D} + B_i}{(-1)^{i+1} 2A_i}$ be the first and second roots of $\Phi_i$, and define $F_i := \frac{(-1)^i}{f_i}$ and $S_i := \frac{(-1)^{i+1}}{s_i}$. Then

$$F_i = \frac{\sqrt{D} + B_i}{2A_{i+1}} \quad \text{and} \quad S_i = \frac{\sqrt{D} - B_i}{2A_{i+1}},$$

with $\partial eg\ F_i > 0 > \partial eg\ S_i$, since $\Phi_i$ are reduced. Furthermore, from the fact that $\Delta_i$ takes $\Phi_i$ to $\Phi_{i+1}$ we know that their roots are related by

$$f_{i+1} = \delta_i - \frac{1}{f_i} \quad \text{and} \quad s_{i+1} = \delta_i - \frac{1}{s_i} \tag{4}$$

Multiplying both by $(-1)^{i+1}$ and using the definition of $F_i$, $S_i$ and $g_i$ we get

$$F_i = g_i + \frac{1}{F_{i+1}} \quad \text{and} \quad S_{i+1} = \frac{1}{g_i + S_i}.$$

Hence

$$F_i = [g_i,\ g_{i+1}, \cdots] \quad \text{and} \quad S_i = [0,\ g_{i-1},\ g_{i-2}, \cdots]$$

Furthermore, using properties of continued fractions we obtain

$$\frac{1}{f_0} = F_0 = [g_0,\ g_1, \ldots,\ g_i,\ F_{i+1}] \tag{5}$$

$$(-1)^{i+1} s_i = \frac{1}{S_i} = \left[g_{i-1},\ g_{i-2}, \cdots,\ g_0,\ \frac{1}{S_0}\right]. \tag{6}$$

**Remark 7.** *Observe that*

$$F_i + S_i = \frac{\sqrt{D}}{A_{i+1}} = [g_i,\ g_{i+1}, \cdots] + [0,\ g_{i-1},\ g_{i-2}, \cdots] \tag{7}$$

**Theorem 6.** *Two properly equivalent reduced indefinite binary quadratic forms belong to the same chain.*

*Proof.* Let $q$ and $Q$ be reduced indefinite binary quadratic forms with coefficients in $\mathbb{Q}((1/T))$ and discriminant $D \neq 0$. Suppose the transformation $H = \begin{pmatrix} \alpha & \beta \\ \gamma & \delta \end{pmatrix} \in SL_2(\mathbb{Q}[T])$ makes them properly equivalent. Since $\alpha\delta - \beta\gamma = 1$, the tuple $(\alpha, \gamma)$ gives a solution to the Diophantine equation $\delta x - \beta y = 1$.Furthermore, in a technical lemma to follow, since $\alpha,\ \beta,\ \delta,\ \gamma$ are entries of a transformation, we have $\partial eg\ \alpha < \partial eg\ \beta$ and $\partial eg\ \beta < \partial eg\ \delta$. Hence $(\alpha, \gamma)$ is the unique non-zero solution in rational polynomials such that $\partial eg\ x < \partial eg\ \beta$ and $\partial eg\ y < \partial eg\ \delta$. Furthermore, from proposition 3, we have that $\gamma/\alpha$ is the $(i-1)^{\text{st}}$ convergent of $\delta/\beta = [g_0,\ g_1, \cdots,\ g_{i-1}]$, if $i$ is even. We can assume $i$ is even, as if it is odd, then the substitution

$$g_{i-2} + \frac{1}{g_{i-1}} = g_{i-2} + \frac{1}{(g_{i-1} - 1) + \frac{1}{1}}$$

extends the continued fraction to even number of terms.
Using the convergents correspondence of the continued fraction of $\delta/\beta$, we get

$$\left[g_0,\ g_1, \cdots,\ g_{i-1}, \frac{1}{F}\right] = \frac{\frac{\delta}{F} + \gamma}{\frac{\beta}{F} + \alpha}.$$

Furthermore, if $F$ and $f$ be the first roots of $Q$ and $q$ respectively, then the transformation $H$, connects them via the identity

$$\frac{1}{f} = \frac{\frac{\delta}{F} + \gamma}{\frac{\beta}{F} + \alpha} = \left[g_0,\ g_1, \cdots,\ g_{i-1}, \frac{1}{F}\right].$$

Since $Q$ is reduced $F$ has a negative degree and $\partial eg\ 1/F > 0$. The continued fraction expansion of $1/f$ is unique up to the $g_{i-1}$ term. On the other hand, from the relations of the roots $f_i$ of the forms in the chain $\Phi_i$, (11), we have $1/f_0 = F_0 = [g_0,\ g_1, \cdots,\ g_{i-1},\ F_i]$. Taking $f_0 = f$, and $F_i = 1/F$, we get $F = 1/F_i = (-1)^i f_i = f_i$, since $i$ is even. Hence $F$ is the first root of the form $\Phi_i$ in the chain where $\Phi_0 = q$.

It remains to show that the second root $s_i$ of $\Phi_i$ is equal to $S$ (the second root of $Q$), given the second root $s_0$ of $\Phi_0$ is equal to $s$ (the second root of $q$). The relations for the second roots $s_i$ of the chain forms given in (12) state

$$(-1)^{i+1} s_i = \frac{1}{S_i} = \left[g_{i-1},\ g_{i-2}, \cdots,\ g_0, \frac{1}{S_0}\right]$$
$$\Rightarrow -s_i = [g_{i-1},\ g_{i-2}, \cdots,\ g_0, -s],$$

since $i$ is even and $s = s_0$. Now, $\partial eg\ s$ is positive, so this expansion is unique up to the $g_0$ term. Furthermore, from proposition 2 applied to the continued fraction of $\delta/\beta$, we know that

$$\frac{\delta}{\gamma} = [g_{i-1},\ g_{i-2}, \cdots,\ g_0] \text{ and } \frac{\beta}{\alpha} = [g_{i-1},\ g_{i-2}, \cdots,\ g_1].$$

Hence from the convergents correspondence, we have

$$-s_i = [g_{i-1},\ g_{i-2}, \cdots,\ g_0, -s] = \frac{-s\delta + \beta}{-s\gamma + \alpha} = -S.$$

The final equality follows from $s$ and $S$ being connected via $H$. Therefore, $S$ is equal to the second root of the form $\Phi_i$ in the chain with $\Phi_0 = q$. Namely, $q$ and $Q$ are in the same chain. □

**Lemma 4.** *If two distinct reduced indefinite binary quadratic forms of the same discriminant $D \neq 0$, are properly equivalent via the transformation $\begin{pmatrix} \alpha & \beta \\ \gamma & \delta \end{pmatrix}$, then*

$$\partial eg\ \alpha \leq \partial eg\ \beta,\ \partial eg\ \gamma < \partial eg\ \delta,\ \text{and } \partial eg\ \beta < \partial eg\ \delta. \tag{8}$$

*Proof.* Since $q$ and $Q$ are properly equivalent, $\alpha\delta = \beta\gamma + 1$. We proceed by case analysis:

Case i. Suppose $\partial eg\ \alpha\delta < 0$. Since $\alpha,\ \delta \in \mathbb{Q}[T]$, $H$ is one of the following

$$\begin{pmatrix} 0 & \pm 1 \\ \mp 1 & \delta \end{pmatrix} \text{ or } \begin{pmatrix} \alpha & \pm 1 \\ \mp 1 & 0 \end{pmatrix}$$

If we are in the latter case, consider $H^{-1} = \begin{pmatrix} 0 & \pm 1 \\ \mp 1 & \alpha \end{pmatrix}$, taking $Q$ to $q$. The matrix $H$ connects the roots by $-\delta = \frac{1}{f} + F$, hence $\partial eg\ \delta = \partial eg\ 1/f > 0$, and since $\partial eg\ \alpha < 0$ and $\partial eg\ \beta = \partial eg\ \gamma = 0$, the conditions are satisfied. For the latter two cases, the conditions are thus satisfied for $H^{-1}$.

If $\beta\gamma = 0$, then $H$ is one of the following

$$\begin{pmatrix} \pm 1 & \beta \\ 0 & \pm 1 \end{pmatrix} \text{ or } \begin{pmatrix} \pm 1 & 0 \\ \gamma & \pm 1 \end{pmatrix}$$

The first transformation connects the first roots $f$ and $F$, by $f - F = \beta$, and since the degrees of both $f$ and $F$ are negative, $\beta = 0$, i.e $H$ is the identity. For the latter matrix, consider the second roots $s$ and $S$. Then $\frac{1}{s} = \gamma + \frac{1}{S}$, and since $s$ and $S$ are of positive degree, we must have $\gamma = 0$, and $H$ is the identity matrix. However, $q \neq Q$, so we can assume that $\beta\gamma \neq 0$.

Case ii. If $\partial eg\ \alpha\delta \geq 0$ and $\beta\gamma \neq \pm 1$. Then $\partial eg\ (\beta\gamma + 1) = \partial eg\ \beta\gamma \geq 0$. Hence

$$\partial eg\ \beta + \partial eg\ \gamma = \partial eg\ \alpha + \partial eg\ \delta \tag{9}$$

and

$$\partial eg\ \alpha < \partial eg\ \beta \Leftrightarrow \partial eg\ \gamma < \partial eg\ \delta \tag{10}$$

(a) Suppose $\partial eg\ \alpha < \partial eg\ \delta$, then (9) implies $\partial eg\ \beta + \partial eg\ \gamma < 2\partial eg\ \delta$.

- if $\partial eg\ \beta = \partial eg\ \gamma$, then $\partial eg\ \beta < \partial eg\ \delta$ and $\partial eg\ \gamma < \partial eg\ \delta$. Then from (10) $\partial eg\ \alpha < \partial eg\ \beta$
- if $\partial eg\ \gamma < \partial eg\ \beta$, then (10) implies that $\partial eg\ \gamma < \partial eg\ \delta$ and $\partial eg\ \alpha < \partial eg\ \beta$. Furthermore, under $H$, the first roots satisfy

  $$\frac{1}{f} = \frac{\gamma + \delta/F}{\alpha + \beta/F}$$

  and since $\partial eg\ f < 0$, we must have $\partial eg\ (\gamma + \delta/F) > \partial eg\ (\alpha + \beta/F)$. Furthermore, $\partial eg\ \gamma < \partial eg\ \delta + \partial eg\ 1/F$, since $\partial eg\ 1/F > 0$. Hence

  $$\partial eg\ (\gamma + \delta/F) = \partial eg\ \delta + \partial eg\ \frac{1}{F} > \partial eg\ (\alpha + \beta/F) \geq \partial eg\ \frac{\beta}{F}$$

  the latter inequality follows from $\partial eg\ \frac{\beta}{F} > \partial eg\ \beta > \partial eg\ \alpha$. Therefore $\partial eg\ \delta > \partial eg\ \beta$.
- if $\partial eg\ \beta < \partial eg\ \gamma$, then (9) implies $\partial eg\ \beta < \partial eg\ \delta$ and $\partial eg\ \alpha < \partial eg\ \gamma$. We use the relation of the second roots under the transformation $H$, namely

  $$\frac{1}{s} = \frac{\gamma + \delta/S}{\alpha + \beta/S}$$
  $$\Rightarrow 1 = \left(\frac{\alpha}{s} - \gamma\right)(\alpha S + \beta)$$

Hence $\partial eg\ \left(\frac{\alpha}{s} - \gamma\right) = -\partial eg\ (\alpha S + \beta)$. Furthermore, $\partial eg\ 1/s < 0$ so

$$\partial eg\ \left(\frac{\alpha}{s} - \gamma\right) = \partial eg\ \gamma = -\partial eg\ (\alpha S + \beta)$$

Furthermore, $\partial eg\ \gamma > \partial eg\ \alpha \geq 0$, i.e $\partial eg\ (\alpha S + \beta) < 0$. Since $\alpha, \beta \in \mathbb{Q}[T]$ and $\partial eg\ S > 0$, this can only happen if $\partial eg\ \alpha S = \partial eg\ \beta$. Hence $\partial eg\ \alpha < \partial eg\ \beta$.

(b) if $\partial eg\ \delta < \partial eg\ \alpha$. Consider $H^{-1}$, taking $Q$ to $q$. Then

$$H^{-1} = \begin{pmatrix} A & B \\ \Gamma & \Delta \end{pmatrix} = \begin{pmatrix} \delta & -\beta \\ -\gamma & \alpha \end{pmatrix}$$

hence $\partial eg\ A < \partial eg\ \Delta$, and the same analysis as in the above cases works.

(c) if $\partial eg\ \alpha = \partial eg\ \delta$, then $2\partial eg\ \alpha = 2\partial eg\ \delta = \partial eg\ \beta + \partial eg\ \gamma$.

- if $\partial eg\ \beta = \partial eg\ \gamma$, then $\partial eg\ \alpha = \partial eg\ \beta = \partial eg\ \gamma = \partial eg\ \delta$. Furthermore, consider

$$1 = \left(\frac{\alpha}{f} - \gamma\right)(\alpha F + \beta) \tag{11}$$

Since, $\partial eg\ 1/f > 0$, and $\partial eg\ F < 0$, we have that

$$-\partial eg\ \beta = -\partial eg\ (\alpha F + \beta) = \partial eg\ \left(\frac{\alpha}{f} - \gamma\right) > \partial eg\ \alpha$$

contradiction.

- if $\partial eg\ \beta > \partial eg\ \gamma$, then $\partial eg\ \alpha < \partial eg\ \beta$ and $\partial eg\ \delta < \partial eg\ \beta$. From (10) we have $\partial eg\ \gamma < \partial eg\ \delta$ and $\partial eg\ \gamma < \partial eg\ \alpha$. Furthermore, taking the degree of (11) we get

$$\partial eg\ (\alpha F + \beta) = -\partial eg\ \left(\frac{\alpha}{f} - \gamma\right)$$

and since $\partial eg\ 1/f > 0$ and $\partial eg\ F < 0$ we have

$$-\partial eg\ \beta = -\partial eg\ (\alpha F + \beta) = \partial eg\ \left(\frac{\alpha}{f} - \gamma\right) > \partial eg\ \alpha$$

But also, $\partial eg\ \beta > \partial eg\ \alpha$, hence $\partial eg\ \alpha < 0$, i.e, $\alpha = 0 = \delta$, and$\beta = \pm 1 = \gamma$, but by assumption $\partial eg\ \beta > \partial eg\ \gamma$, contradiction.

- if $\partial eg\ \beta < \partial eg\ \gamma$, then $\partial eg\ \beta < \partial eg\ \alpha < \partial eg\ \gamma$ and $\partial eg\ \beta < \partial eg\ \delta < \partial eg\ \gamma$. We next consider

$$1 = \left(\frac{\alpha}{s} - \gamma\right)(\alpha S + \beta) \tag{12}$$

Taking degree and using $\partial eg\ 1/s < 0 < \partial eg\ S$, we have

$$\partial eg\ \alpha < \partial eg\ \gamma = -\partial eg\ \left(\frac{\alpha}{s} - \gamma\right) = -\partial eg\ \alpha - \partial eg\ S,$$

i.e $\partial eg\ S < -2\partial eg\ \alpha$ and $\partial eg\ \alpha < 0$. Thus $\alpha = 0$, same analysis as above, gives us a contradiction.

□

# 4 Representation by indefinite binary quadratic forms and the Markov spectrum

**Definition 8.** We say that $A \in \mathbb{Q}((1/T))$ is *represented by an indefinite binary quadratic form* $Q$, if there exist rational polynomials $X, Y$ not both zero, such that $A = Q(X, Y)$.

**Definition 9.** Let $Q$ be an indefinite binary quadratic form of discriminant $D \neq 0$. Let $m(Q) := \inf_{\substack{X,Y \in \mathbb{Q}[T] \\ (X,Y)\neq(0,0)}} \partial eg\ Q(X, Y)$. Then the *Markov spectrum* is defined to be

$$\mathcal{M} := \left\{ \frac{\partial eg\ D(Q)}{2} - m(Q) \mid Q \text{ indefinite binary quadratic form} \right\}.$$

**Proposition 11.** *Properly equivalent binary quadratic forms represent the same elements of* $\mathbb{Q}((1/T))$.

*Proof.* Let $q$ and $Q$ be two binary quadratic forms which are properly equivalent via the transformation $H = \begin{pmatrix} \alpha & \beta \\ \gamma & \delta \end{pmatrix} \in SL(\mathbb{Q}[T])$, and let $M \in \mathbb{Q}((1/T))$ be represented by $q$. That is there are some rational polynomials $x$ and $y$, not both 0, such that $q(x, y) = M$. Let $X = \delta x - \beta y$ and $Y = -\gamma x + \alpha y$, also rational polynomials, then $Q(X, Y) = M$. Finally, $X$ and $Y$, cannot be both zero, since $x$ and $y$ are not both zero and the determinant of $H$ is equal to 1. Therefore $M$ is also represented by $Q$. □

**Theorem 7.** *If the forms* $\left((-1)^i A_i, B_i, (-1)^{i+1} A_{i+1}\right)$*, for $i$ an integer, constitute a chain of reduced forms of discriminant $D \neq 0$, a square in* $\mathbb{Q}((1/T))$*, then the $A_i$'s include all elements of* $\mathbb{Q}((1/T))$ *of degree less than the degree of* $\sqrt{D}$*, which are represented by a form in the chain.*

*Proof.* Let $A \in \mathbb{Q}((1/T))$ with $\partial eg\ A < \partial eg\ \sqrt{D}$ be represented by such a reduced form $Q = \left((-1)^i A_i, B_i, (-1)^{i+1} A_{i+1}\right)$ of discriminant $D$ in a chain. That is there exist rational polynomials $x, y$ not both zero, such that $(-1)^i A_i x^2 + B_i xy + (-1)^{i+1} A_{i+1} y^2 = A$. If we take $\alpha = x$ and $\gamma = y$, where $x, y$ are co-prime, then there exist $\beta, \delta \in \mathbb{Q}[T]$, such that $\alpha\delta - \gamma\beta = 1$. Then the transformation $H = \begin{pmatrix} \alpha & \beta \\ \gamma & \delta \end{pmatrix}$ takes $Q$ to a properly equivalent form $(A, B, C)$ of the same discriminant $D$, which also represents $A$. However, this form is not necessarily reduced. Consider its first and second roots $f = (\sqrt{D} - B)/2A$ and $s = (-\sqrt{D} - B)/2A$. Observe that $\partial eg\ (f - s) = \partial eg\ \sqrt{D} - \partial eg\ A > 0$. Therefore, we can't have both degree of $f$ and $s$ being negative, and we can assume that $\partial eg\ f \geq 0$, otherwise $Q$ is reduced. Furthermore, $\lfloor f \rfloor \neq \lfloor s \rfloor$, so we apply $\begin{pmatrix} 1 & h \\ 0 & 1 \end{pmatrix}$, with $h = \lfloor f \rfloor$. This transformation sends $(A, B, C)$ to $(A, B_1, C_1)$, which is reduced and represents $A$. From theorem 6 $(A, B_1, C_1)$ must be one of the forms in the chain, i.e $A$ must appear amongst the $A_i$'s. □

**Theorem 8.** *Suppose $Q$ is an indefinite binary quadratic form constituting a chain of equivalent reduced forms $\Phi_i = \big((-1)^i A_i, B_i, (-1)^{i+1} A_{i+1}\big)$, for $i \in \mathbb{Z}$. Then*

$$m(Q) = \inf_{i \in \mathbb{Z}} \partial eg\ A_i.$$

*Proof.* Suppose $Q = (A, B, C)$ is a reduced form of discriminant $D$, then

$$\partial eg\ \left(\frac{\sqrt{D}}{A}\right) = \partial eg\ \left(\frac{\sqrt{D} + B}{2A} + \frac{\sqrt{D} - B}{2A}\right) > 0.$$

In particular, $A \in \mathbb{Q}\left((1/T)\right)$ is such that $\partial eg\ A < \partial eg\ \sqrt{D}$. Hence $Q$ represents an element of $\mathbb{Q}\left((1/T)\right)$ of degree smaller than $\partial eg\ \sqrt{D}$. Therefore by theorem 7 is represented by some $A_i$ in the chain of reduced forms equivalent to $Q$. □

# 5 Alternative realisation of the Lagrange and Markov Spectra

In the real case, the Markov and Lagrange Spectra can alternatively be defined via terms of doubly infinite sequences of positive integers. These were first studied by Markov in [5]. He used them to show that the Lagrange spectrum coincides with the Markov spectrum, for numbers below 3. In our setting we work with the analogous object - doubly infinite sequences of polynomials of positive degree, $A = \ldots, g_{-1}, g_0, g_1, \ldots$. To give some intuition on how the Lagrange and Markov spectra is realised via doubly infinite sequences, we reexamine a few identities from sections 2 and 3.

Firstly, suppose we are given $\alpha \in \mathbb{Q}\left((1/T)\right)$, not a rational function, which has a continued fraction expansion $[a_0,\ a_1, \cdots,\ \alpha_{h+1}]$. Then from remark 6 we have the identity

$$\alpha - \frac{p_h}{q_h} = \frac{(-1)^h}{q_h^2\left(\alpha_{h+1} + \frac{q_{h-1}}{q_h}\right)}.$$

Hence $\partial eg\ (\alpha - p_h/q_h) + 2\partial eg\ q_h = -\partial eg\ (\alpha_{h+1} + q_{h-1}/q_h)$. Observe further that

$$\alpha_{h+1} = [a_{h+1},\ a_{h+2}, \cdots] \ \text{ and } \ \frac{q_{h-1}}{q_h} = [0,\ a_h,\ a_{h-1}, \cdots,\ a_1],$$

where each $a_i$ has a positive degree. And the Lagrange constant is given by $l(\alpha) = \limsup_{h\to\infty} \partial eg\ (\alpha_{h+1} + q_{h-1}/q_h)$.

Furthermore, suppose we are given an indefinite binary quadratic form $Q$ of discriminant $D$, which constitutes a chain of equivalent forms

$$\Phi_i = \left((-1)^i A_i, B_i, (-1)^{i+1} A_{i+1}\right).$$

Just as in the discussion after corollary 2, we can define $F_i := \frac{(-1)^i}{f_i}$ and $S_i := \frac{(-1)^{i+1}}{s_i}$, where $f_i$ and $s_i$ are the first and second roots of $\Phi_i$. Then from remark 7 we have

$$F_i + S_i = \frac{\sqrt{D}}{A_{i+1}} = [g_i,\ g_{i+1}, \cdots] + [0,\ g_{i-1},\ g_{i-2}, \cdots],$$

where $g_i$ are rational polynomials of positive degree. Furthermore, the elements of the Markov spectrum are given by $\partial eg\ \sqrt{D} - m(Q)$, which by theorem 8 is the same as $\partial eg\ \sqrt{D} - \inf_{i\in\mathbb{Z}} \partial eg\ A_i$.

Hence it is natural to make the following definition

**Definition 10.** Given a doubly infinite sequence of rational polynomials of positive degree $A = \cdots,\ g_{-1},\ g_0,\ g_1, \cdots$, we define

$$\lambda_i(A) := [g_i,\ g_{i+1}, \cdots] + [0,\ g_{i-1},\ g_{i-2}, \cdots].$$

Furthermore, let

$$L(A) := \limsup_{i\in\mathbb{Z}} \partial eg\ \lambda_i(A) \ \text{ and } \ M(A) := \sup_{i\in\mathbb{Z}} \partial eg\ \lambda_i(A).$$

**Theorem 9.** *The Lagrange spectrum can also be defined as*

$$\mathbb{L} = \{L(A) \mid A \textit{ is a doubly infinite sequence of non-constant rational polynomials}\}.$$

*Proof.* For $\alpha \in \mathbb{Q}\left((1/T)\right)/\mathbb{Q}(T)$, with continued fraction expansion $\alpha = [a_0,\ a_1, \cdots]$, let

$$A = \cdots, a_1,\ a_0,\ a_1, \cdots$$

Then $L(A) = \limsup_{i\to\infty} \partial eg\ \lambda_i(A) = \limsup_{i\to\infty} \deg a_i$. Furthermore, from theorem 1, we know $l(\alpha) = \limsup_{i\to\infty} \deg a_i$. Therefore, $l(\alpha) \in \{L(A) \mid A \text{ as above}\}$.

For the converse, let $A$ be a doubly infinite sequence as in the definition. Then $L(A)$ is either $\limsup_{i\to+\infty} \partial eg\ \lambda_i(A)$ or $\limsup_{i\to-\infty} \partial eg\ \lambda_i(A)$. In the first case we take $\alpha = [g_0,\ g_1, \cdots]$ and in the latter case we take $\alpha = [g_0,\ g_{-1}, \cdots]$. Then $L(A) \in \mathbb{L}$. $\square$

**Theorem 10.** *The Markov spectrum $\mathcal{M}$ can be realised as the set*

$$\mathbb{M} = \{M(A)|A \text{ doubly infinite sequence of non-constant rational polynomials}\}.$$

*Proof.* From the discussion above and remark 7, given an indefinite binary quadratic form $Q$ of discriminant $D \neq 0$ we obtain a doubly infinite sequence of non constant polynomials $A$, such that $M(A) = \sqrt{D} - m(Q)$. Hence $\mathcal{M} \subseteq \mathbb{M}$.

On the other hand, given a doubly infinite sequence of rational polynomials of positive degree $A = \cdots,\ g_{-1},\ g_0,\ g_1, \cdots$, we consider

$$\lambda_i(A) = [g_i,\ g_{i+1}, \cdots] + [0,\ g_{i-1},\ g_{i-2}, \cdots] \in \mathbb{Q}\left((1/T)\right).$$

Thus we can find an element of $\mathbb{Q}\left((1/T)\right)$, say $A_{i+1}$, of degree $-\partial eg\ g_i < 0$, such that $\lambda_i(A) = 1/A_{i+1}$. Let $F_i = [g_i,\ g_{i+1}, \cdots]$ and $S_i = [0,\ g_{i-1},\ g_{i-2}, \cdots]$, then $F_i + S_i = 1/A_{i+1}$. Choose $B_i \in \mathbb{Q}\left((1/T)\right)$, such that

$$F_i = \frac{1 + B_i}{2A_{i+1}} \quad \text{and} \quad S_i = \frac{1 - B_i}{2A_{i+1}}.$$

Then we consider $f_i = (-1)^i/F_i$ and $s_i = (-1)^i/S_i$, i.e

$$f_i = \frac{1 - B_i}{2(-1)^i a_i} \quad \text{and} \quad s_i = \frac{1 + B_i}{2(-1)^i a_i},$$

where $4A_{i+1}a_i = 1 - B_i^2$. Furthremore, $\partial eg\ S_i < 0 < \partial eg\ F_i$ and thus $\partial eg\ f_i < 0 < \partial eg\ s_i$. Therefore, $f_i$ and $s_i$ are the roots of the reduced indefinite binary quadratic form $Q_i = \left((-1)^i a_i, B_i, (-1)^{i+1}A_{i+1}\right)$ of discriminant 1. From the continued fraction expansion of $F_i$ and $S_i$ we have

$$F_i = g_i + \frac{1}{F_{i+1}} \quad \text{and} \quad \frac{1}{S_i} = g_{i-1} + S_{i-1}$$
$$\Rightarrow f_{i+1} = \delta_i - \frac{1}{f_i} \quad \text{and} \quad s_{i+1} = \delta_i - \frac{1}{s_i},$$

where $\delta_i = (-1)^i g_i$. Then the transformation $\Delta_i = \begin{pmatrix} 0 & 1 \\ -1 & \delta_i \end{pmatrix}$ sends $Q_i$ to $Q_{i+1}$, and in particular $a_{i+1} = A_{i+1}$. Hence the forms $Q_i = \left((-1)^i A_i, B_i, (-1)^{i+1}A_{i+1}\right)$ are reduced, of discriminant 1 and in a chain. From theorem 8 we know $\inf_{i\in\mathbb{Z}} A_i = m(Q)$, where $Q$ is indefinite quadratic form of discriminant 1 properly equivalent to $Q_i$. Then

$$\begin{aligned} M(A) = \sup_{i\in\mathbb{Z}} \partial eg\ \lambda_i(A) &= \sup_{i\in\mathbb{Z}} \partial eg\ \left(\frac{1}{A_{i+1}}\right) \\ &= -\inf_{i\in\mathbb{Z}} \partial eg\ A_{i+1} \\ &= -m(Q). \end{aligned}$$

Hence $\mathbb{M} \subseteq \mathcal{M}$. □

**Corollary 3.** *The Markov Spectrum* $\mathcal{M} = \mathbb{N} \cup \{\infty\}$.

*Proof.* From the above theorem

$\mathcal{M} = \mathbb{M} = \{M(A) | A \text{ doubly periodic seq of non-constant rational polynomials}\}$.

Furthermore, $M(A) = \sup_{i \in \mathbb{Z}} \partial eg\ \lambda_i(A)$, and $\partial eg\ \lambda_i(A) = \partial eg\ g_i$, where $g_i \in \mathbb{Q}[T]$ has positive degree. The result follows. □

**Corollary 4.** *The Lagrange and Markov spectra are the same.*